\documentclass[letterpaper, 10 pt, conference]{ieeeconf}

\IEEEoverridecommandlockouts

\usepackage{amsmath,amssymb,amsfonts}
\usepackage{algorithmic}
\usepackage{cite}
\usepackage[utf8]{inputenc}
\usepackage[T1]{fontenc}
\usepackage{graphicx}
\usepackage{textcomp}
\usepackage{amsthm}
\usepackage{xcolor}
\usepackage{color}
\usepackage{subcaption}
\usepackage{dsfont}
\usepackage{soul}
\usepackage{mathrsfs}

\usepackage{multirow}

\usepackage[most]{tcolorbox}
\newtcolorbox{highlighted}{colback=yellow,coltext=red,breakable}

\newcommand{\mR}{{\mathbb R}}

\newcommand{\mP}{{\mathbb P}}
\newcommand{\mK}{{\mathbb K}}

\newcommand{\bR}{{\mathbf R}}

\newcommand{\indicator}{\mathds{1}}

\newcommand{\bs}{{\mathbf s}}

\newcommand{\intzeroinfty}{\int_0^{\infty}}
\usepackage{float}
\usepackage{amsmath}
\usepackage{subcaption}
\usepackage{array}

\usepackage{color}
\newcommand{\yu}{\color{black}}

\usepackage{multirow}

\begin{document}

\title{Data-driven optimal control under safety constraints using sparse Koopman approximation}

\author{Hongzhe Yu, Joseph Moyalan, Umesh Vaidya, and Yongxin Chen
\thanks{Financial support from NSF under grants 1942523, 2008513, 2031573 and NSF CPS award 1932458 is greatly acknowledged.}
\thanks{H. Yu and Y.\ Chen are with the School of Aerospace Engineering,
Georgia Institute of Technology, Atlanta, GA; {\{hyu419,yongchen\}@gatech.edu}}
\thanks{J. Moyalan and U. Vaidya are with the Department of Mechanical Engineering, Clemson University, Clemson, SC; {\{jmoyala,uvaidya\}@clemson.edu}}
}

\maketitle
\thispagestyle{empty}
\pagestyle{empty}

\begin{abstract}
In this work we approach the dual optimal reach-safe control problem using sparse approximations of Koopman operator. Matrix approximation of Koopman operator needs to solve a least-squares (LS) problem in the lifted function space, which is computationally intractable for fine discretizations and high dimensions. The state transitional physical meaning of the Koopman operator leads to a sparse LS problem in this space. Leveraging this sparsity, we propose an efficient method to solve the sparse LS problem where we reduce the problem dimension dramatically by formulating the problem using only the non-zero elements in the approximation matrix with known sparsity pattern. The obtained matrix approximation of the operators is then used in a dual optimal reach-safe problem formulation where a linear program with sparse linear constraints naturally appears. We validate our proposed method on various dynamical systems and show that the computation time for operator approximation is greatly reduced with high precision in the solutions. 

\end{abstract}

\section{Introduction}\label{introduction}
{\yu In the era of machine learning and data-driven technology, we often face problems where models for complex systems are difficult to obtain while abundant data collected from the system are available, such as biology science\cite{yeung2015global}, finance, social networks, and fluid dynamics\cite{brunton2016discovering}.} Data-driven system identification \cite{champion2019data, brunton2016discovering, quade2018sparse} of complex dynamical systems has thus seen huge advancements in the recent years. {\yu For nonlinear dynamics the most commonly used method to identify a system is by solving a least-squares regression in the span of nonlinear basis functions.} Linear operator theory revolving around Koopman and Perron-Frobenius (PF) operators \cite{lasota1998chaos} is a powerful tool for analyzing nonlinear system in the \textit{`lifted'} {\yu nonlinear function} spaces, where the operator \textit{`lifted'} dynamics becomes linear and describes individual or collective movement. The operators can be approximated using finite dimensional linear mappings (matrices) in a least-square sense using data. The state-of-the-art method for approximating Koopman and PF operators in this manner are the Dynamic Mode Decomposition (DMD) and its extensions \cite{mezic2005spectral, tu2013dynamic, williams2015data, proctor2016dynamic, huang2018data}.

{\yu However, the above mentioned data-driven approximations and identifications for nonlinear dynamics involves a linear regression under the hood,} which means the computation complexity is intractable with increasingly higher dimensions. Many methods in the literature seek to compute the linear finite dimension operator approximation efficiently. In \cite{sinha2020computationally} the authors proposed a method based on Cholesky decomposition to reduce the dimension of the matrix being inverted. Other methods include nonlinear model reduction \cite{alla2017nonlinear} and exploring dynamics sparsity structures to decouple the system into different interconnected subsystems \cite{schlosser2022sparsity}.

In this work we start from the physical interpretation of the Koopman operator to approach the scalability issue. Koopman operator describes the linear system evolution in the lifted space, which corresponds to a sparse linear state transition matrix approximation with a known sparsity pattern for pre-defined discretization grids. Starting from this observation, we formulate the sparse least-squares problem using only the known nonzero elements in the sparse matrix, which reduced the problem dimension by a large factor. The problem then becomes an equivalent linear system of equations of much lower dimension than the original one. The sparse approximations are then used in a dual optimal control formulation \cite{yu2022data, prajna2004nonlinear}, serving as a sparse linear constraint in a linear program for optimal control synthesis.
Different from existing {\yu sparse system identification method such as SINDy \cite{brunton2016discovering}}, the proposed method {\yu starts} directly from a known sparse pattern arising from the operators' physical meaning {\yu instead of using $L_1$ norm to promote sparsity, and is not restricted to polynomial dynamics.} The proposed method does not need system reduction or matrix manipulations {\yu when} solving the least-squares problem. 

The rest of the paper is structured as follows. Section II provides a brief introduction to our framework’s necessary
preliminaries. Section III talks about sparse approximations of the linear operators. The construction and solving of sparse least squre problem is explained in Section IV. In Section V, we present some simulation results followed by {\yu the} conclusion in Section VI.

\section{background and notations}\label{chp:background}
In this section we briefly introduce Koopman and PF operators and their finite dimensional approximations.
\subsection{Koopman and Perron-Frobenius Operator}
\label{sec:Koopman}
Consider the dynamical system
\begin{eqnarray}
\dot x(t)=f(x(t)),\;\; x(t)\in X, \;\; x(0) = x_0.\label{sys}
\label{dyn_sys}
\end{eqnarray}
We use $s_t(x_0)$ and $x(t)$ interchangeably to denote the mapping from initial state $x_0$ to the solution of system (\ref{dyn_sys}) at time $t$. $\bs_{-t}(x)$ represents the set $\bs_{-t}(x) = \{y: s_t(y)=x\}$. ${\mathds{1}}_{A}(x)$ denotes the indicator function on a set $A$. Koopman and PF operators are tools to study \eqref{dyn_sys} in lifted spaces \cite{korda2018convergence}. 

Koopman operator $\mathbb{K}_t $ for system~\eqref{sys} is defined as 
\begin{eqnarray}[\mathbb{K}_t \varphi](x)=\varphi(s_t(x)), \label{eq:defn_koopman_operator}
\end{eqnarray}
where $\varphi$ is a test function. The infinitesimal generator for Koopman operator $\mathbb{K}_t$ is defined as
\begin{eqnarray}
\lim_{t\to 0}\frac{[\mathbb{K}_t\varphi](x)-\varphi(x)}{t}=f(x)\cdot \nabla \varphi(x)=:{\cal K}_f \varphi. 
\label{eq:infinitesimal_Koopman}
\end{eqnarray}
PF operator $\mathbb{P}_t$ for \eqref{dyn_sys} is defined as 
\vspace{-0.08in}
\begin{equation}
    \int_{\bs_{-t}(A)} \psi(x) dx = \int_{A} \mP_t[\psi](x) dx,~~ \forall A\subset X,
\end{equation}
and its infinitesimal generator is given by 
\begin{equation}
\lim_{t\to 0}\frac{[\mathbb{P}_t\psi](x)-\psi(x)}{t}=-\nabla \cdot (f(x) \psi(x))=: {\cal P}_f\psi. \label{eq:PF_generator}
\end{equation}
The duality between the two operators reads
\begin{equation}\label{eq:duality}
\int_{X}[\mathbb{K}_t \varphi]\psi dx=
\int_{X}[\mathbb{P}_t \psi]\varphi dx,~ \forall \psi, \varphi.
\end{equation}
\subsection{Extended Dynamic Mode Decomposition (EDMD)}
Numerical methods were proposed to approximate Koopman operator $\mK_t$ in a finite dimensional setting. Extendede Dynamic Mode Decomposition (EDMD) \cite{williams2015data} approximates $\mathbb{K}_t$ within a linearly-spanned space $\mathcal{D}$ of nonlinear basis
\begin{equation}
    \Psi(x) \triangleq [\psi_1(x),\dots,\psi_N(x)]^T.
\end{equation}
It seeks a matrix representation $K$ of $\mathbb{K}_t$ with respect to $\Psi$ by minimizing
\vspace{-0.08in}
\begin{equation}
\left \lVert \Psi_{y} -  K\Psi_{x} \right \rVert_F^2,
\label{eq:EDMD}
\end{equation}
where $\Psi_{x}=[\Psi(x_{1}) \ldots \Psi(x_{M})], \Psi_{y}=[\Psi(y_{1}) \ldots \Psi(y_{M})]$ are the lifted data in dimension $N \times M$, $M$ being the number of data points. Commonly used basis include polynomial \cite{yu2021convex} and Gaussian {\yu radial basis function (RBF)}. The minimizer to problem \eqref{eq:EDMD} {\yu is}
\begin{equation}
    K^{\star}\! =\! (\frac{1}{M}\!\sum_{k=1}^{M}\Psi(y_k)\Psi(x_k)^T) (\frac{1}{M}\!\sum_{k=1}^{M} \Psi(x_k)\Psi(x_k)^T)^{\dagger}.
    \label{eq:pseudo_inv}
\end{equation}
Computing the pseudo inverse in \eqref{eq:pseudo_inv} becomes intractable for large size matrices. In this work we explore the structure of the problem to reduce the computation complexity.

\subsection{Dual optimal reach-safe control formulation}\label{sec:ocp problem}
The operator approximations lead to a dual formulation to optimal control. We consider control-affine dynamics
\begin{equation}
     \dot x = f + g u,
     \label{eq:control-affine}
\end{equation}
where $u(x) \in \mR^m$ is the feedback policy. For a fixed $u(x)$, denote $s_t(x_0)$ or $x(t)$ the solution to \eqref{eq:control-affine} at time $t$. We seek an optimal policy $u(x)$ to drive the system from set $X_0$ to set $X_r$ while avoiding (unsafe) set $X_u$, formulated as
\begin{subequations}\label{eq:optimal_control_problem}
\begin{align}
    \inf_{u(\cdot)} &\;\;\int_{X} \intzeroinfty l(\bs_t(x_0), u(\bs_t(x_0))) dt h_0(x_0)dx_0\\
{\rm s.t.} &\;\;\;\;\intzeroinfty \indicator_{X_u} (\bs_t(x_0)) dt = 0,~\forall x_0 \in X_0. 
\end{align}
\label{eq:optimal_control_problem}
\end{subequations}
Here {\yu $l$ is the running cost}, $h_0$ denotes the initial distribution of $x_0.$
It is showed in \cite{yu2022data} that, for system \eqref{eq:control-affine}, \eqref{eq:optimal_control_problem} can be reformulated into
\begin{equation}
\begin{split}
\inf_{\rho, u} &\;\;\int_{X}(q(x) + \lVert u(x)\rVert_1) \rho(x)dx \\ 
    \textrm{s.t.} &\;\;\;\; \nabla\cdot [({f + g u})\rho](x)=h_0(x)\\
    &\;\;\;\;\int_{X} \indicator_{X_u}(x) \rho(x)dx = 0,
\end{split}
\label{eq:l1-problem}
\end{equation}
where $q(x)$ is the state cost and $\lVert \cdot \rVert_1$ denotes the $1-$norm.
Here a new variable termed \textit{`occupation measure'} is introduced and represented by $\rho(x)$. It is defined as $\rho(x) \triangleq \intzeroinfty [\mP_t h_0](x)dt$.
Problem \eqref{eq:l1-problem} is bi-linear in $u$ and $\rho$. By introducing a new variable $\bar\rho \triangleq \rho u$, \eqref{eq:l1-problem} is turned into a convex problem
\begin{equation}\label{eq:LP_convex}
\begin{split}
\inf_{\rho, \bar\rho} &\;\;\int_{X} q(x)\rho(x) + \lVert \bar\rho(x)\rVert_1 dx \\ 
    \textrm{s.t.} &\;\;\;\; \nabla\cdot ({f \rho + g \bar\rho})(x)=h_0(x)\\
    &\;\;\;\;\int_{X} \indicator_{X_u}(x) \rho(x)dx = 0
\end{split}
\end{equation}
in variable{\yu s} $(\rho, \bar\rho)$. After solving \eqref{eq:LP_convex}, $u$ can be recovered using $u(x)=\frac{\bar\rho(x)}{\rho(x)}$. We will show that \eqref{eq:LP_convex} will become a linear program using operator approximations.

\section{Sparse approximations of operators}
\label{sec:sparse_approx}
 We investigate the finite {\yu dimensional} approximation of the Koopman operator and  explore the sparsity pattern in this approximation arising from its physical meaning. 
\subsection{EDMD for controlled dynamical systems}
To clarify the notations, in \eqref{eq:control-affine} we denote $g = [g_1,\dots,g_m]\in\mR^{n\times m}$ and $u=[u_1,\dots,u_m]^T \in \mR^{m}$. Data ($\Psi_{x}, \Psi_{y}$) are collected from simulated system trajectories.
For simplicity, we collect a set of $m+1$ data, $\{(\Psi_{x}^k, \Psi_{y}^k)\}_{k=0}^m$ which contains $1$ data from autonomous system ($u=0$), and $m$ data from system with different inputs, the $k^{th}$ data corresponding to the one-hot input $e_k$ system. We solve \eqref{eq:EDMD} for each $k$ to get $K_0, \dots, K_m$. Note that these can also be approximated jointly \cite{yu2022data}. For PF operator approximations, for a function $\phi$ {\yu with coefficients $C_{\phi}$ defined as $\phi\triangleq\Psi^T C_{\phi}$,} we approximate \cite{das2018data}
\begin{equation}
    [\mathbb{K}_t \phi] \approx \Psi^T (K^T  C_{\phi}),\;\; [\mathbb{P}_t \phi] \approx \Psi^T (P^T  C_{\phi}),
    \label{eq:approx_PF}
\end{equation}
and by definition \eqref{eq:PF_generator},
\begin{equation}
    \mathcal{P}_{f} \approx \frac{P^T-I}{\Delta t}.
    \label{eq:approx_PF_generator}
\end{equation}
Define the integral $\Lambda \triangleq \int \Psi(x)\Psi(x)^T dx$, and in view of the duality \eqref{eq:duality}, $\forall y_1\triangleq \Psi^T C_{y_1}, y_2\triangleq\Psi^TC_{y_2}$, we have 
\begin{eqnarray*}
    \langle [\mathbb{K}_ty_1], y_2 \rangle 
    &=& C_{y_1}^T K \Lambda C_{y_2}\\
    &=&\langle y_1, [\mathbb{P}_t y_2] \rangle \\
    &=& C_{y_1}^T \Lambda P^T  C_{y_2}.
\label{eq:duality_Psi}
\end{eqnarray*}
The above is true for all $y_1, y_2$, {\yu so we have} 
\begin{equation}
    P^T  = \Lambda^{-1}K\Lambda. 
    \label{eq:relation_K_P}
\end{equation}
Note that this is true for all $K_i$ and $P_i$. When the basis functions $\Psi$ are orthogonal, we assume that $\Lambda$ is diagonal dominant and $P^T=\Lambda^{-1}K\Lambda \approx K.$

Now we approximate the constraints in (\ref{eq:LP_convex}). we parameterize the variables $\rho(x) \triangleq  \Psi(x)^Tv$, and $ \bar\rho(x) \triangleq [w_1, \dots, w_m]^T\Psi(x)$ within ${\cal{D}}$
and approximate
\begin{equation}
\begin{split}
    \nabla\cdot ({f \rho + g \bar\rho}) 
    &=\nabla\cdot (f \rho) + \sum_{i=1}^{m}\nabla \cdot (g_i\rho u_i)\\
    &= -\mathcal{P}_{f}\rho - \sum_{i=1}^m \mathcal{P}_{g_i}\bar\rho_i\\
    &\approx \Psi^T (\frac{I - P _0^T}{\Delta t} v) + \Psi^T(\frac{I - (P _i-P _0)^T}{\Delta t} w_i).
\end{split}
\label{eq:approx_constraint}
\end{equation}

\label{sec:sparsity_structure}


\subsection{Sparse least-squares problem}
\begin{figure}
\begin{subfigure}{.24\textwidth}
  \centering
  \includegraphics[width=.95\linewidth]{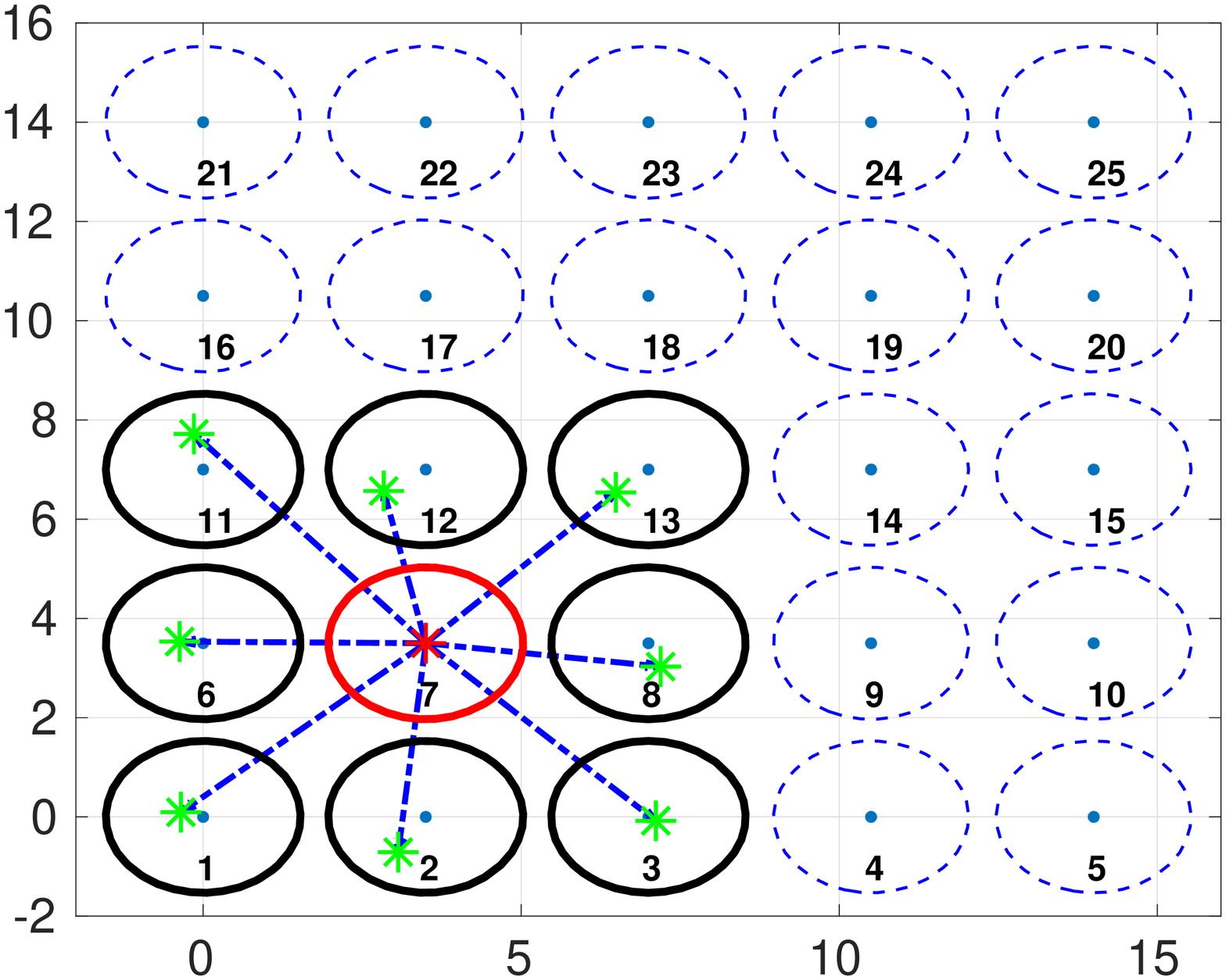}
  \caption{Grids and state transitions.}
  \label{fig:sys_trans_sfig1}
\end{subfigure}
\begin{subfigure}{.24\textwidth}
  \centering
  \includegraphics[width=.95\linewidth]{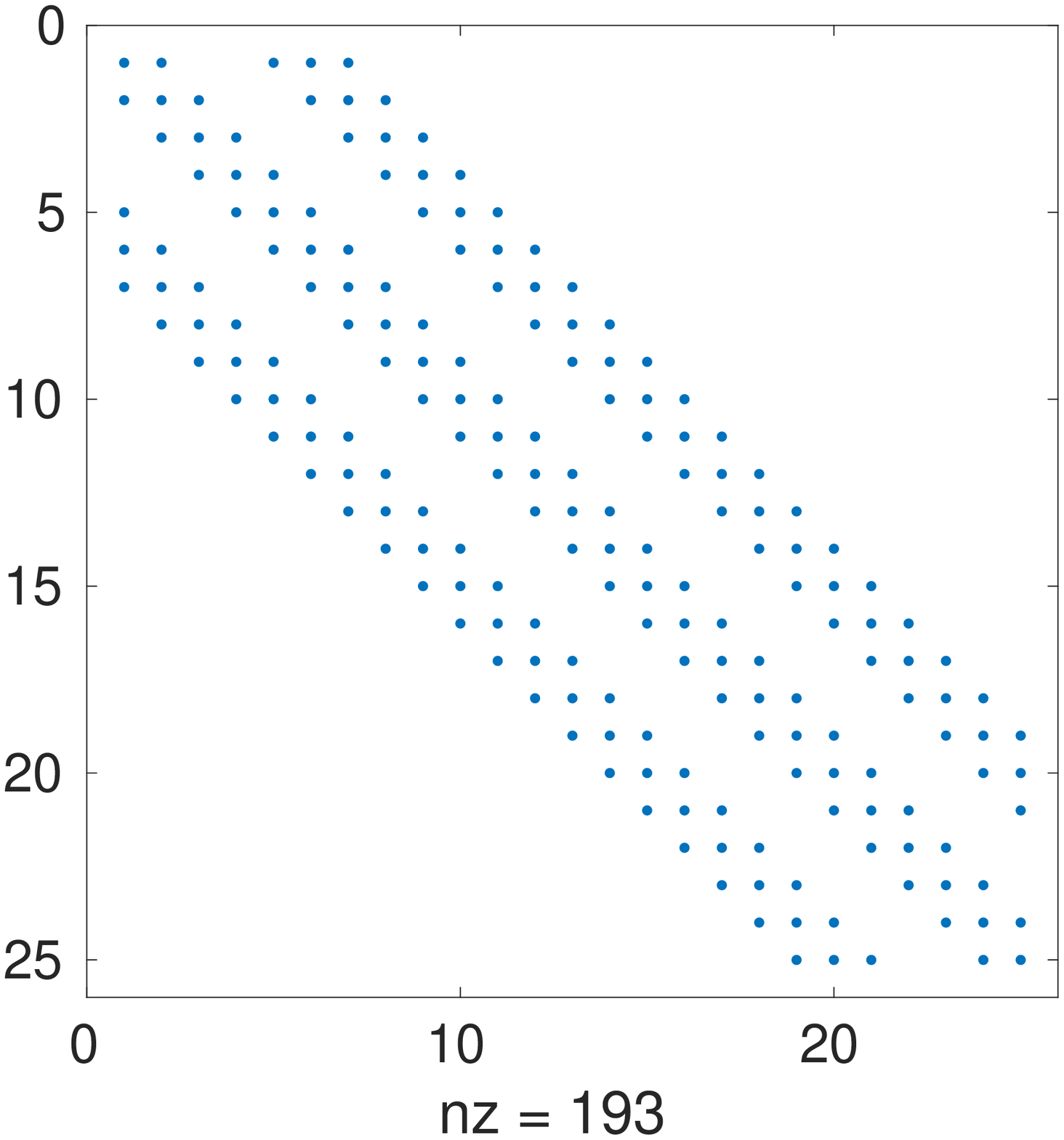}
  \caption{Sparse construction of $K$.}
  \label{fig:sys_trans_sfig2}
\end{subfigure}
\caption{Grid layout and sparsity. {\yu Fig.\ref{fig:sys_trans_sfig1} shows the basis functions layout, and Fig.\ref{fig:sys_trans_sfig2} shows the sparsity pattern of the resulting Koopman matrix.}}
\label{fig:sys_trans}
\end{figure}

\begin{figure}[h]
    \centering
    \includegraphics[width=0.45\textwidth]{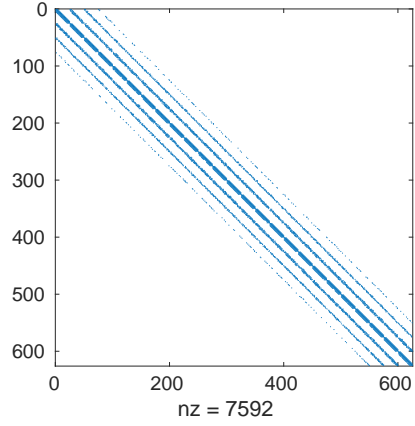}
    \caption{A typical sparsity pattern (after truncating elements below $1e^{-6}$) of $K$ by solving \eqref{eq:EDMD} using \eqref{eq:pseudo_inv}. It has 7592 non-zero entries in a 625$\times$625 matrix. Sampling time $\Delta t=10^{-3}$.}
    \label{fig:K1_truncated}
\end{figure}

Both Koopman and PF operators describe the system evolution. Specifically, $\mathbb{K}$ captures the evolution of a point in the state space, and $\mathbb{P}$ captures the evolution of the distribution of a collection of states. The matrix approximation $K$ in \eqref{eq:EDMD} has a sparse structure due to this physical meaning, since within a small sampling time the range of system evolution is limited. We discretize the state space as shown in Fig. \ref{fig:sys_trans_sfig1} and use Gaussian RBF basis, then $K$ will have a banded-diagonal sparsity structure representing neighbouring-grid activations as shown in Fig. \ref{fig:sys_trans_sfig2}. Fig. \ref{fig:K1_truncated} is a typical truncated $K$ obtained from the pseudo-inverse solution \eqref{eq:pseudo_inv}.

With known banded-diagonal sparsity pattern, we now formulate a sparse LS problem from \eqref{eq:EDMD} using only the non-zero elements in $K$. Assuming that $K$ has $z$ non-zero entries which are collected in a vector $v \in \bR^{z}$. For notation brevity, we write the operator which constructs $K$ from $v$ as $\mathscr{K}(\cdot)$ or $\mathscr{K}_{\cdot}$. Its adjoint operator is denoted as $\mathscr{K}^{+}(\cdot)$ or $\mathscr{K}^{+}_{\cdot}$ and is defined by

\begin{equation}
    \langle \mathscr{K}_{v}, L \rangle = \langle v, \mathscr{K}_{L}^+ \rangle, \forall{L}
    \label{eq:adj_operatorK}
\end{equation}
{\yu for any matrix} $L$ {\yu which} has the same dimension as $\mathscr{K}_{v}$. Since
\begin{equation}
    \begin{split}
        \langle \mathscr{K}_{v}, L \rangle &=\sum_{i,j, K_{i,j}\neq 0} {K_{i,j}L_{i,j}}\\
        &=\sum_{i=1}^{z} (v_{K})_i (v_{L})_i =\langle v, v_{L} \rangle,
    \end{split}
\end{equation}
we know immediately that $\mathscr{K}_{L}^+ = v_{L}, \forall L$,
where $v_{L} \in \mR^z$ is constructed from $L$ using the same rule for $v$ from $K$. For instance, for a diagonal matrix $K$, $\mathscr{K}_{K}^{+} = diag(K)$ and $\mathscr{K}_{v} = diag(v)$,
where the function $diag(\cdot)$ takes either matrix or vector inputs, and extracts the diagonal elements to a vector in the former case and construct a diagonal matrix from the input vector in the latter case. We use similar mappings for the banded diagonal structure in Fig. \ref{fig:K1_truncated}. Writing problem \eqref{eq:EDMD} in $v$ as
\begin{equation}
    \begin{split}
        J_{E}(v)&=\lVert \mathscr{K}_{v}\Psi_x - \Psi_y \rVert_F^2\\
        &= \langle \mathscr{K}_{v}\Psi_x, \mathscr{K}_{v}\Psi_x \rangle - 2\langle \mathscr{K}_{v}\Psi_x, \Psi_y \rangle + \langle \Psi_y, \Psi_y\rangle\\
        &= {\rm Tr}[\mathscr{K}_{v}^T\mathscr{K}_{v}\Psi_x\Psi_x^T] - 2{\rm Tr}[\mathscr{K}_{v}\Psi_x\Psi_y^T],
    \end{split}\nonumber
    \label{eq:sparse_obj}
\end{equation}
where from the second to the third equation we omit the terms irrelevant to $v$. Now taking the first order approximation of $J_E(v)$ with respect to {\yu a} perturbation $\delta v$
\begin{equation}
    \begin{split}
        &\phantom{{}\leq{}}J_{E}(v + \delta v)\\
        &= {\rm Tr}[\mathscr{K}_{v + \delta v}^T\mathscr{K}_{v + \delta v}\Psi_x\Psi_x^T] - 2{\rm Tr}[\mathscr{K}_{v + \delta v}\Psi_x\Psi_y^T]\\
        &={\rm Tr}[(\mathscr{K}_{v}^T + \mathscr{K}_{\delta v}^T)(\mathscr{K}_{v} + \mathscr{K}_{\delta v})\Psi_x\Psi_x^T] -\\
        &\;\;\;\;\;2{\rm Tr}[(\mathscr{K}_{v}+\mathscr{K}_{\delta v})\Psi_x\Psi_y^T]\\
        &\approx J_E(v)+\langle \mathscr{K}_{\delta v}, 2(\mathscr{K}_{v}\Psi_x\Psi_x^T-\Psi_y\Psi_x^T)\rangle\\
        &=J_E(v) + \langle \delta v, 2(\mathscr{K}^+(\mathscr{K}_{v}\Psi_x\Psi_x^T)-\mathscr{K}^+(\Psi_y\Psi_x^T)) \rangle\\
        &=J_E(v) + \langle \delta v, \frac{\partial J_E}{\partial v} \rangle,
    \end{split}
\end{equation}
where we used the linearity of operator $\mathscr{K}$, Taylor's first-order approximation, and definition \eqref{eq:adj_operatorK}.
{\yu From the last line we know that} the minimizer $v^*$ solves the equation 
\begin{equation}
    \mathscr{K}^+(\mathscr{K}_{v}\Psi_x\Psi_x^T)-\mathscr{K}^+(\Psi_y\Psi_x^T) = 0.
\label{eq:sp_normeq}
\end{equation}
Note the matrix $\mathscr{K}_{v} \in \mR^{N\times N}$. Equation \eqref{eq:sp_normeq} is not explicit in $v$, but the LHS is a linear function of $v$. We thus seek to solve an equivalent linear equation in $v$
\begin{equation}
    S v = \mathscr{K}^+(\Psi_y\Psi_x^T)
    \label{eq:sp_normeq_lin}
\end{equation}
by first constructing $S \triangleq [S_1, \dots, S_z].$ The $i$th column $S_i$ is constructed by assigning $v=e_i$ in the LHS of \eqref{eq:sp_normeq}
\begin{equation}
    S_i = \mathscr{K}^+(\mathscr{K}_{e_i}\Psi_x\Psi_x^T).
    \label{eq:S_columns}
\end{equation}
{\yu Now we} transformed the least-squares problem \eqref{eq:EDMD} of dimension $N \times N$ into a linear equation system \eqref{eq:sp_normeq_lin} of dimension $z$. $K$ is then obtained by $K = \mathscr{K}_v$.
The PF operator and its generator can then be approximated using \eqref{eq:relation_K_P} and \eqref{eq:approx_PF_generator}. 

\section{Constructing and solving the sparse least-squares problem}\label{sec:Algo}
Section \ref{sec:sparse_approx} constructed a sparse least-squares problem for operator approximation. In this section we propose an algorithm to construct and solve problem \eqref{eq:sp_normeq_lin}.
\subsection{The {\yu size} of full and sparse problems}
Following the grid layout shown in Fig. \ref{fig:sys_trans_sfig1}, the sparse matrix construction will have a banded-diagonal structure as shown in Fig. \ref{fig:sys_trans_sfig2}. {\yu Assume} the layout of the basis functions has $n_b$ basis along each dimension, then the exact non-zero diagonals are fixed,  assuming the system can only go to its neighbouring grid within the sampling time. {\yu For instance, }in 2 dimensional case, $K$ has $3\times 3$ nonzero diagonals corresponding to the 9 bold basis shown in Fig. \ref{fig:sys_trans_sfig1}. Similarly, for 3 dimensional system, $K$ has $3\times 3^2$ nonzero diagonals. The number of nonzero elements for a $n$ dimensional system is less than $3^n\times (n_b)^n$, compared with the size of the matrix $(n_b)^{2n}$. The ratio between the two decreases exponentially with $n$. For banded-diagonal {\yu matrices,} the mapping $\mathscr{K}^+$ and $\mathscr{K}$ are both known. $\mathscr{K}^+$ is the concatenation of the nonzero diagonals into a vector, and $\mathscr{K}$ is reversing operation on the input vector to construct the bands. 
\subsection{Constructing the sparse linear equation}
\paragraph{Sparse data lifting}
In both the full matrix setting \eqref{eq:EDMD} and the sparse setting \eqref{eq:sp_normeq_lin}, lifted data matrices $(\Psi_x, \Psi_y)$ are computed where we evaluate the basis functions $\Psi=[\Psi_1 \dots \Psi_N]^T$ at each data point $x_{k}$ and $y_{k}, k=1,\dots M$. The dimension of $\Psi_{x}$ and $\Psi_{y}$ grows exponentially with the state dimension, which makes the lifting step alone a computation bottle neck. One observation is that the activated $\Psi$ is also sparse. Each column of $\Psi_{x}$ and $\Psi_{y}$ is sparse as the basis functions are near zero at the positions far from the data, assuming the basis are almost orthogonal to each other as in the case of RBF. With this, we construct a sparse $\Psi_{x}^s$ and $\Psi_{y}^s$ by selecting only a fixed number of closest basis to a given data point. The terms $\Psi_{y}^s{\Psi_{x}^s}^T$ and $\Psi_{x}^s{\Psi_{x}^s}^T$ in \eqref{eq:sp_normeq_lin} are also sparse after this sparse lifting. Often $\Psi_{y}^s{\Psi_{x}^s}^T$ and $\Psi_{x}^s{\Psi_{x}^s}^T$ become singular, {\yu with which the pseudo-inverse method} \eqref{eq:EDMD} is not viable while the proposed method can leverage the sparse $\Psi^s$. 

\paragraph{Sparse least-squares problem}
The problem of interest is \eqref{eq:sp_normeq_lin}. The solving takes only a small portion of the time and constructing \eqref{eq:sp_normeq_lin} from \eqref{eq:sp_normeq} needs careful design to reduce the overall time consumed compared with full pseudo inverse method. We seek to construct $S$ using $\eqref{eq:S_columns}$. The key to an efficient construction of $S$ is the observation that the matrix $\mathscr{K}_{e_i}$ only has $1$ nonzero element, assuming at position $(l,s)$. The multiplication $\mathscr{K}_{e_i}\Psi_x\Psi_x^T$ selects the $s^{th}$ row of $\Psi_x\Psi_x^T$ as the $l^{th}$ row in the resulting matrix. The position $(l, s)$ is efficiently found using sorting and binary search. The operator $\mathscr{K}^+$ then turns the one row matrix back into the vector $S_i$. As we have already seen, $\Psi_x\Psi_x^T$ is sparse by construction, which means that the matrix $S$ will also be sparse, and \eqref{eq:sp_normeq_lin} is a sparse linear system of equations.
After constructing \eqref{eq:sp_normeq_lin}, then iterative methods such as conjugate gradient descent \cite{shewchuk1994introduction} is used to solve this sparse linear equation.

\section{Experiments}\label{section:experiments}
In this section experiments and analysis are presented. We first introduce an LP formulation for the dual optimal control problem. We then compare the computation time for $K$ between the full pseudo-inverse and the proposed sparse method. Finally we use the obtained sparse $K$ to solve the optimal reach-safe control problem. 
\subsection{$L_1$ input-regularized problem and linear programming}
\label{sec:LP_SDP}
After solving the sparse least-squares problem (\ref{eq:sp_normeq_lin}) for $v$ and using operator $\mathscr{K}^+$, we obtain sparse matrix approximations $K$ and $P$ of Koopman and PF operator, respectively. We then formulate the problem \eqref{eq:LP_convex} into a linear program (LP) with sparse constraints. By definition \eqref{eq:PF_generator}, and with a direct replacement of \eqref{eq:approx_constraint} into \eqref{eq:LP_convex}, we write \eqref{eq:LP_convex} as an LP 
\begin{equation}
\begin{split}
\inf_{C_{\rho}, w_1,\dots,w_m, s_1,\dots,s_m} &\;\;I_q^TC_{\rho} + \sum_{i=1}^mI_{\Psi}^Ts_i + I_m^TC_{\rho} \\ 
    \textrm{s.t.} \;\;\frac{I - P_0^T}{\Delta t}C_{\rho} + &\; \sum_{i=1}^m  \frac{I - (P_i-P_0)^T}{\Delta t}w_i = C_h\\
    &I_{X_u}^TC_{\rho} = 0\\
    &s_i \geq \lvert w_i \rvert, i=1, \dots, m\\
    &\lvert w_i \rvert \leq u_{\rm max} C_{\rho}, i=1, \dots, m
\end{split}
\label{eq:LP}
\end{equation}
where $\rho=C_{\rho}^T\Psi, u_i=w_i^T\Psi, h=C_h^T\Psi$, and
\begin{equation}
    \begin{split}
        &I_q = \int_{X} q\cdot\Psi , \;\; I_{\Psi} = \int_{X} \Psi , \;\; I_{X_u} = \int_{X} \Psi\cdot\indicator_{X_u}
    \end{split}
\end{equation}
are constants, $s_i$ are slack variables, and $P_i$ are sparse matrix approximations of PF generators \eqref{eq:approx_PF_generator}. Notice that in \eqref{eq:LP} we introduce an additional off-road cost term $I_m$. {\yu When an off-road cost map is available, such as shown in Fig. \ref{fig:2d_nav}, $I_m$ is the cost on the map grids.} We also add additional constraints $\lvert w_i \rvert \leq u_{\textrm{max}} C_{\rho}$ which is an enforcement of bounded control inputs $\lvert u \rvert \leq u_{\rm max}$ corresponding to actuation limitations. This constraint guarantees \eqref{eq:LP} can have non-trivial solutions, because without it the minimum value is zero by choosing $C_\rho=0, w_i=0, s_i=0, \forall{i}$.

\subsection{Sparse Koopman matrix approximation}
\paragraph{Computation time}
We compare the time consumed for solving \eqref{eq:EDMD} for different number of basis functions on a 2D single-integrator system
\begin{equation}
    \dot{x} = u
    \label{eq:single_integrator}
\end{equation}
for $x=[x_1,x_2]^T, u=[u_1,u_2]^T,$ and for the 2D Van der Pol 
\begin{equation}
    \dot{x}_1 = x_2; \;\; \dot{x}_2 = -x_1+x_2(1-x_1^2) + u.
    \label{eq:vdp}
\end{equation}
Solving time for full and sparse method and Frobenius norm of the difference between the two are concluded in table \ref{tab:2d_compare_LS}. 
\begin{table}[h]
    \begin{subtable}[h]{0.45\textwidth}
        \centering
        \begin{tabular}{|c|c|c|c|c|c|}
\hline
\multicolumn{1}{|c|}{\multirow{2}{*}{\# basis}} 
 &\multicolumn{2}{c|}{\# nonzero elements}                         & \multicolumn{2}{c|}{solving time (s)}                         & \multicolumn{1}{c|}{\multirow{2}{*}{norm diff.}} \\ \cline{2-5}
\multicolumn{1}{|c|}{}                  & \multicolumn{1}{c|}{full} & \multicolumn{1}{c|}{sparse} & \multicolumn{1}{c|}{full} & \multicolumn{1}{c|}{sparse} & \multicolumn{1}{c|}{}                  \\ \hline
\multicolumn{1}{|c|}{25}                  & \multicolumn{1}{c|}{390625} & \multicolumn{1}{c|}{5329} & \multicolumn{1}{c|}{9.86} & \multicolumn{1}{c|}{\textbf{0.91}} & \multicolumn{1}{c|}{$3.03e^{-8}$}                  \\\hline
\multicolumn{1}{|c|}{30}                  & \multicolumn{1}{c|}{810000} & \multicolumn{1}{c|}{7744} & \multicolumn{1}{c|}{40.23} & \multicolumn{1}{c|}{\textbf{2.21}} & \multicolumn{1}{c|}{$1.04e^{-8}$}                  \\ \hline
\multicolumn{1}{|c|}{40}                  & \multicolumn{1}{c|}{2560000} & \multicolumn{1}{c|}{13924} & \multicolumn{1}{c|}{117.82} & \multicolumn{1}{c|}{\textbf{3.39}} & \multicolumn{1}{c|}{$1.88e^{-8}$}                  \\ \hline  
\multicolumn{1}{|c|}{60}                  & \multicolumn{1}{c|}{12960000} & \multicolumn{1}{c|}{31684} & \multicolumn{1}{c|}{765.81} & \multicolumn{1}{c|}{\textbf{14.69}} & \multicolumn{1}{c|}{$7.41e^{-8}$}            \\
\hline
\end{tabular}
       \caption{2D Integrator dynamics \eqref{eq:single_integrator}}
       \label{tab:3d_compute_integrator}
    \end{subtable}
    \hfill
\begin{subtable}[h]{0.45\textwidth}
    \centering
    \begin{tabular}{|c|c|c|c|c|c|}
\hline
\multicolumn{1}{|c|}{\multirow{2}{*}{\# basis}} 
 &\multicolumn{2}{c|}{\# nonzero elements}                         & \multicolumn{2}{c|}{solving time (s)}                         & \multicolumn{1}{c|}{\multirow{2}{*}{norm diff.}} \\ \cline{2-5}
\multicolumn{1}{|c|}{}                  & \multicolumn{1}{c|}{full} & \multicolumn{1}{c|}{sparse} & \multicolumn{1}{c|}{full} & \multicolumn{1}{c|}{sparse} & \multicolumn{1}{c|}{}                  \\ \hline
\multicolumn{1}{|c|}{25}                  & \multicolumn{1}{c|}{390625} & \multicolumn{1}{c|}{43557} & \multicolumn{1}{c|}{1.95} & \multicolumn{1}{c|}{\textbf{0.8}} & \multicolumn{1}{c|}{$2.12e^{-8}$}                  \\\hline
\multicolumn{1}{|c|}{30}                  & \multicolumn{1}{c|}{810000} & \multicolumn{1}{c|}{62500} & \multicolumn{1}{c|}{47.14} & \multicolumn{1}{c|}{\textbf{3.4}} & \multicolumn{1}{c|}{$2.09e^{-8}$}                  \\ \hline
\multicolumn{1}{|c|}{40}                  & \multicolumn{1}{c|}{2560000} & \multicolumn{1}{c|}{115600} & \multicolumn{1}{c|}{131.3} & \multicolumn{1}{c|}{\textbf{5.67}} & \multicolumn{1}{c|}{$1.09e^{-8}$}                  \\ \hline  
\multicolumn{1}{|c|}{60}                  & \multicolumn{1}{c|}{12960000} & \multicolumn{1}{c|}{270400} & \multicolumn{1}{c|}{854.7} & \multicolumn{1}{c|}{\textbf{21.65}} & \multicolumn{1}{c|}{$4.78e^{-8}$}                \\
\hline
\end{tabular}
    \caption{Van der Pol dynamics \eqref{eq:vdp}}
    \label{tab:3d_compute_drift}
 \end{subtable}
 \caption{Comparing for Koopman approximations for two different systems. Solving time includes data lifting, constructing and solving \eqref{eq:sp_normeq_lin}. $40000$ data are used. Results are averaged from 5 experiments. {\yu \textit{`norm diff.'} represents the Frobenius norm of the difference of the matrices obtained from two methods.}}
\label{tab:2d_compare_LS}
\end{table}

\paragraph{Computation time and precision for 3D system}
For higher dimensional systems, the full pseudo-inverse {\yu ground truth} took too long to {\yu obtain the results while} the proposed sparse solution can be computed efficiently with high precision. We {\yu show results of} the 3D single integrator system \eqref{eq:single_integrator} for $x=[x_1, x_2, x_3]^T, u=[u_1,u_2,u_3]^T$. Because the lack of {\yu ground truth, we compute the prediction error $\lVert K\Psi_{x} - \Psi_{y} \rVert_F^2$ instead of the Frobenius norm of the difference from ground truth.}

After solving the LS problem \eqref{eq:EDMD} for a $K$, we compute the error $\lVert K\Psi_{x} - \Psi_{y} \rVert_F^2$. Table \ref{tab:3d_compute_integrator} summarizes the results.
\begin{table}[h]
        \centering
        \begin{tabular}{|c|c|c|c|c|c|}
\hline
\# basis & size & nnz & solve time (s) & $\Delta t$ (s) & pred. err \\ \hline
$20$ & $6.4e^{7}$  & $195059$ & $21.75$ & $3.2e^{-7}$ &  $7.27e^{-8}$ \\ \hline
$25$ & $2.4e^{8}$ & $389017$ & $51.82$  & $2.4e^{-8}$ &  $6.78e^{-10}$ \\ \hline
$30$ & $7.29e^{8}$  &  $681493$ & $574.92$ & $4e^{-10}$ & $3.13e^{-13}$ \\  \hline
\end{tabular}
 \caption{Computation time and prediction error for 3D single integrator \eqref{eq:single_integrator}. Number of data is fixed at $8e^4$ for all experiments. Prediction error is measured by $\lVert K \Psi_{x} - \Psi_{y} \rVert_F^2$.}
\label{tab:3d_compute_integrator}
\end{table}

\subsection{Reach-safe optimal control problem}
\begin{figure}[h]
\begin{subfigure}{.49\textwidth}
    \centering
    \includegraphics[width=.95\textwidth]{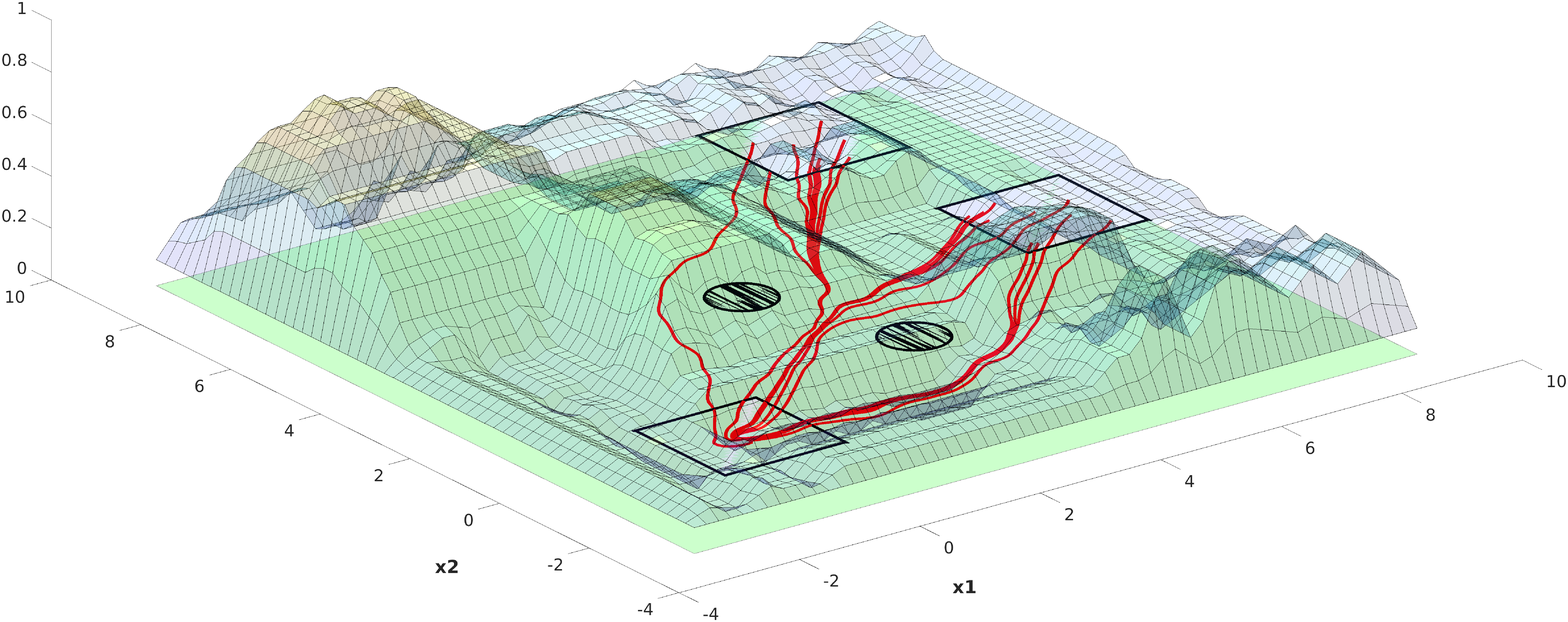}
    \caption{{\yu Cost map, initial states, target states,} and optimized trajectories.}
    \label{fig:2d_nav}
    \end{subfigure}
\begin{subfigure}{.24\textwidth}
  \centering
  \includegraphics[width=.94\linewidth]{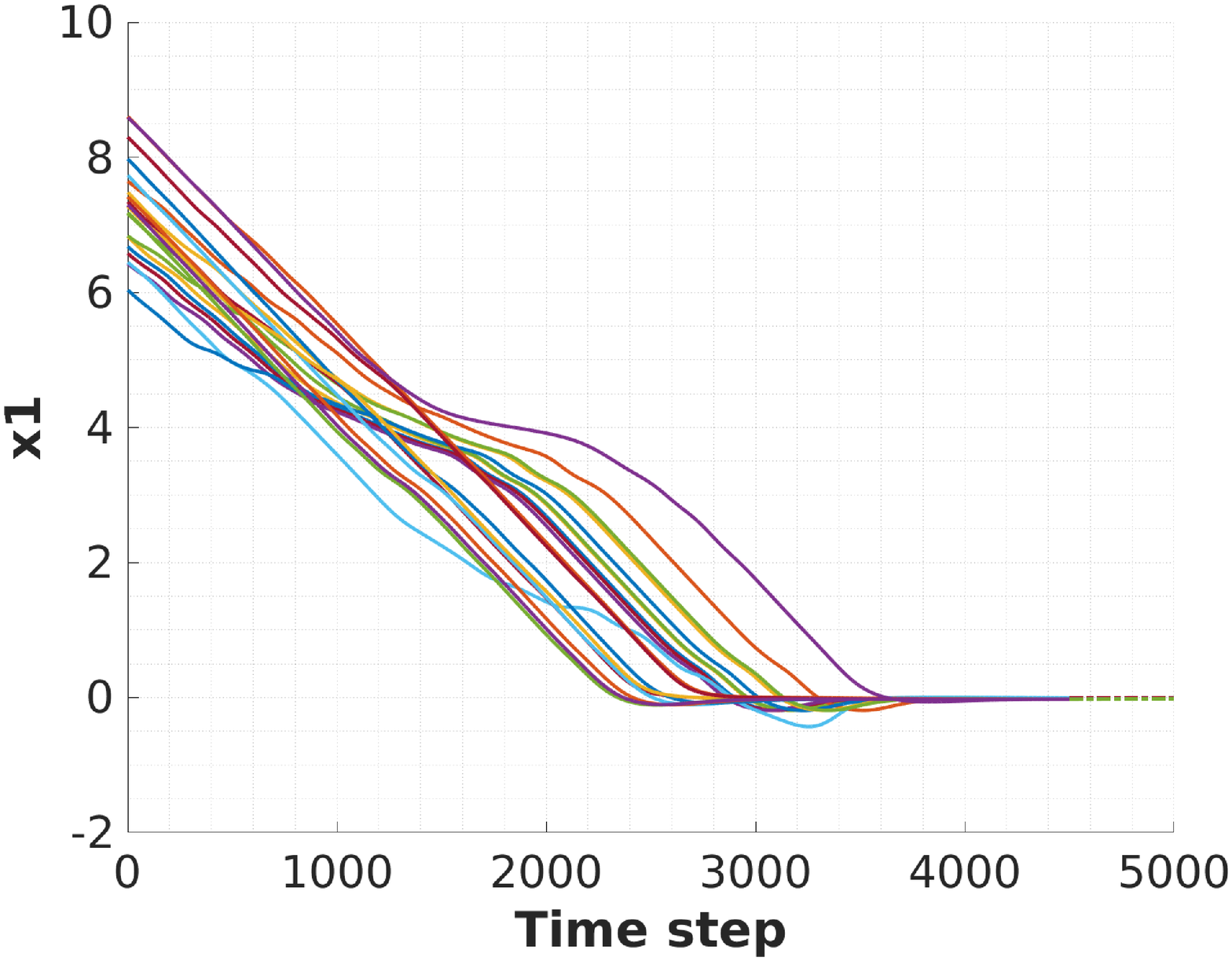}
  \caption{ $x_1$ of sampled trajectories}
  \label{fig:2d_inte_x1}
\end{subfigure}
\begin{subfigure}{.24\textwidth}
  \centering
  \includegraphics[width=.94\linewidth]{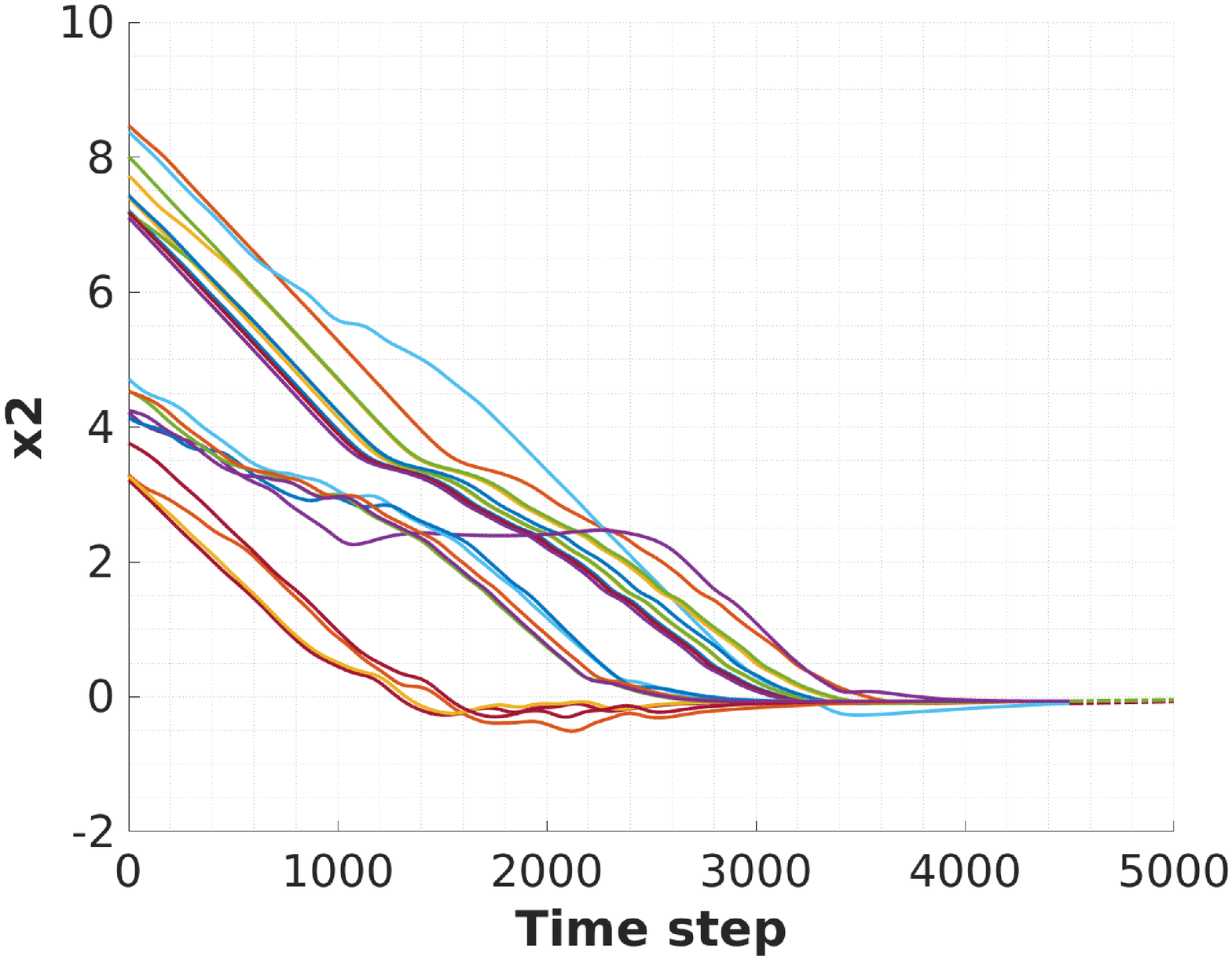}
  \caption{ $x_2$ of sampled trajectories}
  \label{fig:2d_inte_x2}
\end{subfigure}
\caption{Navigation problem with pre-defined cost map in problem \eqref{eq:LP} for a 2d system.}
\end{figure}

We use the sparse $K$ and $P$ to solve problem \eqref{eq:LP}.
\paragraph{2D off-road navigation}
We consider the problem \eqref{eq:LP} for the 2D single integrator dynamics \eqref{eq:single_integrator} with a pre-defined cost map $I_m$ in the following setting
\begin{itemize}
    \item $X_0^1$: A box $[6,8]\times[7,9]$; $X_0^2$: A box $[7,9]\times[3,5]$
    \item $X_u^1$: A circle of radius 1, centered at $(3,4)$
    \item $X_u^2$: A circle of radius 1, centered at $(4, 1.5)$
    \item $X_r$: A circle of radius 1, centered at $(0,0)$.
\end{itemize}
Fig. \ref{fig:2d_nav} shows the optimized system trajectories, and Fig. \ref{fig:2d_inte_x1} - Fig. \ref{fig:2d_inte_x2} show the state plots.

\paragraph{Van der Pol Oscillator} Consider the controlled Van der Pol dynamics \eqref{eq:vdp}. The experiment settings are
\begin{itemize}
    \item $X_0$: A box $[-5,5]\times[-5,5]$
    \item $X_u^1$: A circle of radius 0.5 and center $(-1.5,1)$
    \item $X_u^2$: A circle of radius 1, centered at $(1.5, -1)$
        \item $X_u^3$: A circle of radius 0.5, centered at $(-1, -3)$
    \item $X_r$: A circle of radius 0.5, centered at $(0,0)$.
\end{itemize}

\begin{figure}[ht]
  \centering
  \begin{subfigure}[b]{0.45\textwidth}
  \centering
  \includegraphics[width=0.99\linewidth]{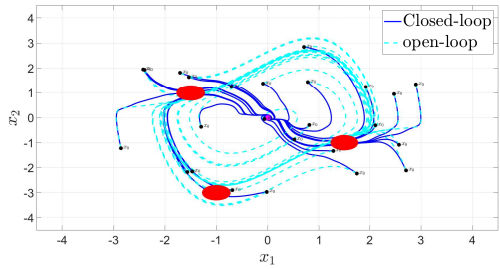}
\caption{Trajectories reaching target without collision constraints.}
\label{no_obstacle}
\end{subfigure}
\begin{subfigure}[b]{0.45\textwidth}
  \centering
  \includegraphics[width=0.9\linewidth]{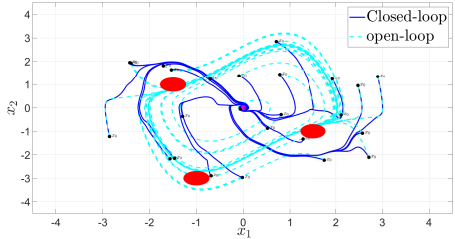}
\caption{Reach-safe trajectories.}
\label{obstacle}
\end{subfigure}
\label{obs}
\caption{Optimal reach-safe results for Van der Pol dynamics. Light blue lines are open-loop trajectories, and dark blue lines are the optimal control resulted trajectories.}
\end{figure}

\begin{figure}[ht]
  \centering
  \begin{subfigure}[b]{0.45\textwidth}
  \centering
  \includegraphics[width=0.9\linewidth]{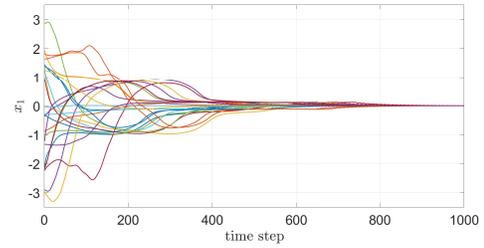}
\caption{state $x_1$}
\label{x_1}
\end{subfigure}
\begin{subfigure}[b]{0.45\textwidth}
  \centering
  \includegraphics[width=0.9\linewidth]{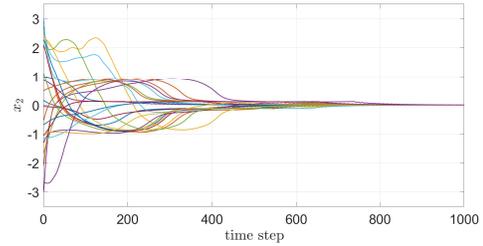}
\caption{state $x_2$}
\label{x_2}
\end{subfigure}
\label{xvdp_states}
\caption{States plot for Van der Pol oscillator reach avoid problem corresponding to Fig.~\ref{obstacle}.}
\end{figure}

\begin{figure}[h]
    \centering
    \includegraphics[width=.8\linewidth]{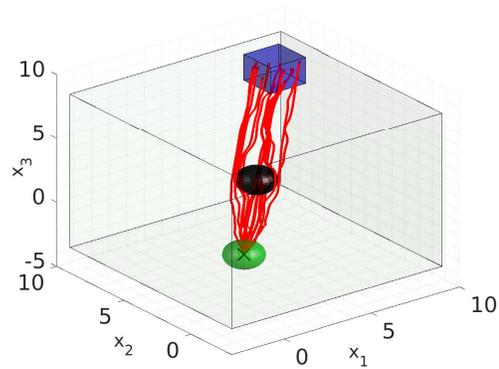}
    \caption{3D-integrator reach-safe problem. Blue box is the set $X_0$, black circle indicates the unsafe region $X_u$, and the green circle is the target set $X_r$. }
    \label{fig:3d_inte}

\end{figure}

\paragraph{3D single integrator}
Consider the simple 3D integrator dynamics \eqref{eq:single_integrator}. The experiment settings are
\begin{itemize}
    \item Initial distributed set: A box $[6.5,8.5]^3$
    \item Unsafe region: A circle of radius 1, centered at $(3,3,3)$
    \item Target region: A circle of radius 1, centered at $(0,0,0)$.
\end{itemize}
After the trajectories enter the target region, we switch to a local LQR controller. Fig. \ref{fig:3d_inte} shows the optimized trajectories.

\section{Conclusion}
In this work we explore the sparsity in the Koopman operator approximations which arises from the state transitional physical meaning of Koopman operator. Using only the nonzero elements in the sparse matrix with known sparsity, we transform the least-squares problem into a sparse linear system of equations. The obtained operator approximations are then used in a dual optimal control  formulation which leads to a linear programming with sparse linear constraints. Results show that our sparse method is much more efficient than the pseudo-inverse method while preserving high precision in the solution.

\bibliographystyle{IEEEtran}
\bibliography{reference}

\begin{thebibliography}{10}
\providecommand{\url}[1]{#1}
\csname url@samestyle\endcsname
\providecommand{\newblock}{\relax}
\providecommand{\bibinfo}[2]{#2}
\providecommand{\BIBentrySTDinterwordspacing}{\spaceskip=0pt\relax}
\providecommand{\BIBentryALTinterwordstretchfactor}{4}
\providecommand{\BIBentryALTinterwordspacing}{\spaceskip=\fontdimen2\font plus
\BIBentryALTinterwordstretchfactor\fontdimen3\font minus
  \fontdimen4\font\relax}
\providecommand{\BIBforeignlanguage}[2]{{%
\expandafter\ifx\csname l@#1\endcsname\relax
\typeout{** WARNING: IEEEtran.bst: No hyphenation pattern has been}%
\typeout{** loaded for the language `#1'. Using the pattern for}%
\typeout{** the default language instead.}%
\else
\language=\csname l@#1\endcsname
\fi
#2}}
\providecommand{\BIBdecl}{\relax}
\BIBdecl

\bibitem{yeung2015global}
E.~Yeung, J.~Kim, J.~Gon{\c{c}}alves, and R.~M. Murray, ``Global network
  identification from reconstructed dynamical structure subnetworks:
  Applications to biochemical reaction networks,'' in \emph{2015 54th IEEE
  Conference on Decision and Control (CDC)}.\hskip 1em plus 0.5em minus
  0.4em\relax IEEE, 2015, pp. 881--888.

\bibitem{brunton2016discovering}
S.~L. Brunton, J.~L. Proctor, and J.~N. Kutz, ``Discovering governing equations
  from data by sparse identification of nonlinear dynamical systems,''
  \emph{Proceedings of the national academy of sciences}, vol. 113, no.~15, pp.
  3932--3937, 2016.

\bibitem{champion2019data}
K.~Champion, B.~Lusch, J.~N. Kutz, and S.~L. Brunton, ``Data-driven discovery
  of coordinates and governing equations,'' \emph{Proceedings of the National
  Academy of Sciences}, vol. 116, no.~45, pp. 22\,445--22\,451, 2019.

\bibitem{quade2018sparse}
M.~Quade, M.~Abel, J.~Nathan~Kutz, and S.~L. Brunton, ``Sparse identification
  of nonlinear dynamics for rapid model recovery,'' \emph{Chaos: An
  Interdisciplinary Journal of Nonlinear Science}, vol.~28, no.~6, p. 063116,
  2018.

\bibitem{lasota1998chaos}
A.~Lasota and M.~C. Mackey, \emph{Chaos, fractals, and noise: stochastic
  aspects of dynamics}.\hskip 1em plus 0.5em minus 0.4em\relax Springer Science
  \& Business Media, 1998, vol.~97.

\bibitem{mezic2005spectral}
I.~Mezi{\'c}, ``Spectral properties of dynamical systems, model reduction and
  decompositions,'' \emph{Nonlinear Dynamics}, vol.~41, no.~1, pp. 309--325,
  2005.

\bibitem{tu2013dynamic}
J.~H. Tu, ``Dynamic mode decomposition: Theory and applications,'' Ph.D.
  dissertation, Princeton University, 2013.

\bibitem{williams2015data}
M.~O. Williams, I.~G. Kevrekidis, and C.~W. Rowley, ``A data--driven
  approximation of the koopman operator: Extending dynamic mode
  decomposition,'' \emph{Journal of Nonlinear Science}, vol.~25, no.~6, pp.
  1307--1346, 2015.

\bibitem{proctor2016dynamic}
J.~L. Proctor, S.~L. Brunton, and J.~N. Kutz, ``Dynamic mode decomposition with
  control,'' \emph{SIAM Journal on Applied Dynamical Systems}, vol.~15, no.~1,
  pp. 142--161, 2016.

\bibitem{huang2018data}
B.~Huang and U.~Vaidya, ``Data-driven approximation of transfer operators:
  Naturally structured dynamic mode decomposition,'' in \emph{2018 Annual
  American Control Conference (ACC)}.\hskip 1em plus 0.5em minus 0.4em\relax
  IEEE, 2018, pp. 5659--5664.

\bibitem{sinha2020computationally}
S.~Sinha, S.~P. Nandanoori, and E.~Yeung, ``Computationally efficient learning
  of large scale dynamical systems: A koopman theoretic approach,'' in
  \emph{2020 IEEE International Conference on Communications, Control, and
  Computing Technologies for Smart Grids (SmartGridComm)}.\hskip 1em plus 0.5em
  minus 0.4em\relax IEEE, 2020, pp. 1--6.

\bibitem{alla2017nonlinear}
A.~Alla and J.~N. Kutz, ``Nonlinear model order reduction via dynamic mode
  decomposition,'' \emph{SIAM Journal on Scientific Computing}, vol.~39, no.~5,
  pp. B778--B796, 2017.

\bibitem{schlosser2022sparsity}
C.~Schlosser and M.~Korda, ``Sparsity structures for koopman and
  perron--frobenius operators,'' \emph{SIAM Journal on Applied Dynamical
  Systems}, vol.~21, no.~3, pp. 2187--2214, 2022.

\bibitem{yu2022data}
H.~Yu, J.~Moyalan, U.~Vaidya, and Y.~Chen, ``Data-driven optimal control of
  nonlinear dynamics under safety constraints,'' \emph{IEEE Control Systems
  Letters}, vol.~6, pp. 2240--2245, 2022.

\bibitem{prajna2004nonlinear}
S.~Prajna, P.~A. Parrilo, and A.~Rantzer, ``Nonlinear control synthesis by
  convex optimization,'' \emph{IEEE Transactions on Automatic Control},
  vol.~49, no.~2, pp. 310--314, 2004.

\bibitem{korda2018convergence}
M.~Korda and I.~Mezi{\'c}, ``On convergence of extended dynamic mode
  decomposition to the {K}oopman operator,'' \emph{Journal of Nonlinear
  Science}, vol.~28, no.~2, pp. 687--710, 2018.

\bibitem{yu2021convex}
H.~Yu, J.~Moyalan, D.~Tellez-Castro, U.~Vaidya, and Y.~Chen, ``Convex optimal
  control synthesis under safety constraints,'' in \emph{2021 60th IEEE
  Conference on Decision and Control (CDC)}.\hskip 1em plus 0.5em minus
  0.4em\relax IEEE, 2021, pp. 4615--4621.

\bibitem{das2018data}
A.~K. Das, B.~Huang, and U.~Vaidya, ``Data-driven optimal control using
  perron-frobenius operator,'' \emph{arXiv preprint arXiv:1806.03649}, 2018.

\bibitem{shewchuk1994introduction}
J.~R. Shewchuk \emph{et~al.}, ``An introduction to the conjugate gradient
  method without the agonizing pain,'' 1994.

\end{thebibliography}

\end{document}